\documentclass[final,3p,times]{elsarticle}

\usepackage{amssymb}
\usepackage{amsmath}
\usepackage{colortbl}
\usepackage{subcaption}
\usepackage{float}
\usepackage{amsmath}
\usepackage{amsfonts}
\usepackage{xcolor}
\usepackage{graphicx}
\usepackage{soul}

\journal{}

\begin{document}

\begin{frontmatter}



\title{A car-following framework for traffic instability and lane changes}


\author{Nicholas Mankowski}
\author{Hassan Mushtaq}
\author{Hanliang Guo\corref{cor1}}
\ead{hguo@owu.edu}

\affiliation{organization={Department of Mathematics and Computer Science\\ Ohio Wesleyan University},
            city={Delaware},
            state={OH},
            postcode={43015},
            country={USA}}

\begin{abstract}
This paper develops a computational framework based on a car-following model to study traffic instability and lane changes.
Building upon Newell's classical first-order car-following model, we show that, both analytically and numerically, there exists a vehicle-density-dependent critical reaction time that determines the stability of single-lane traffic.
Specifically, perturbations to the equilibrium system decay with time for low reaction time and grow for high reaction time. This critical reaction time converges to Newell's original result in the continuum limit.
Additionally, we propose a psychology-based lane-changing mechanism that builds a quantitative connection between the driver's psychological factor (frustration level) and the driving condition. 
We show that our stochastic lane-changing model can faithfully reproduce interesting phenomena like load-balancing of different lanes. 
Our model supports the result that more frequent lane changes only marginally benefit the driver's overall velocity.
\end{abstract}




\date{}

\end{frontmatter}



\section{Introduction}

The traffic dynamics is a complex system problem in which everyone’s driving decisions affect the overall traffic flow. Research on traffic dynamics began as early as the 1930s, with early contributions such as Greenshields’ study on traffic flow characteristics~\cite{greenshields1935study}, and a surge in the number of publications on the topic was seen since the late 1950s, driven by foundational works like Lighthill and Whitham’s and Richards' kinematic wave theory of traffic flow~\cite{lighthill1955kinematic,richards1956shock} and Chandler's pioneer work in car-following models~\cite{chandler1958traffic}.
Models of traffic dynamics not only help us understand these emergent behaviors but also enable the development of intelligent traffic management systems~\cite{helbing2001traffic}, autonomous vehicle algorithms~\cite{treiber2013traffic}, and sustainable urban planning solutions~\cite{levinson2007planning}. By leveraging advanced tools like high-resolution simulation~\cite{nagel1992cellular}, researchers and engineers are unlocking insights that can make transportation safer, more efficient, and less environmentally impactful. The interplay between human psychology, engineering, and complexity science makes traffic modeling a rich field for discovery, as demonstrated by studies like those on driver behavior modeling~\cite{toledo2007driving} and congestion propagation in urban networks~\cite{saberi2020simple,ambuhl2023understanding}.
Loosely speaking, the models of traffic dynamics can be divided into macroscopic models and microscopic models. Macroscopic models describe the traffic flow in terms of aggregate variables, such as density, flow rate, and average velocity. They are often based on the hydrodynamic theory of traffic flow, which views traffic as a fluid. The most significant example of macroscopic models is perhaps the Lighthill-Whitham-Richards (LWR) model~\cite{lighthill1955kinematic,richards1956shock}, which treats traffic as a compressible fluid. The parameters of LWR model are obtained from the fundamental diagram, determining the flow rate at any point on the road by the traffic density at that point. While macroscopic models provide high computational efficiency and analytical traceability via their concise representation of traffic flow, they are limited in the microscopic detail and, perhaps more importantly, cannot capture the stochastic nature of real-life traffic dynamics.

\begin{figure}
    \includegraphics[width=\linewidth]{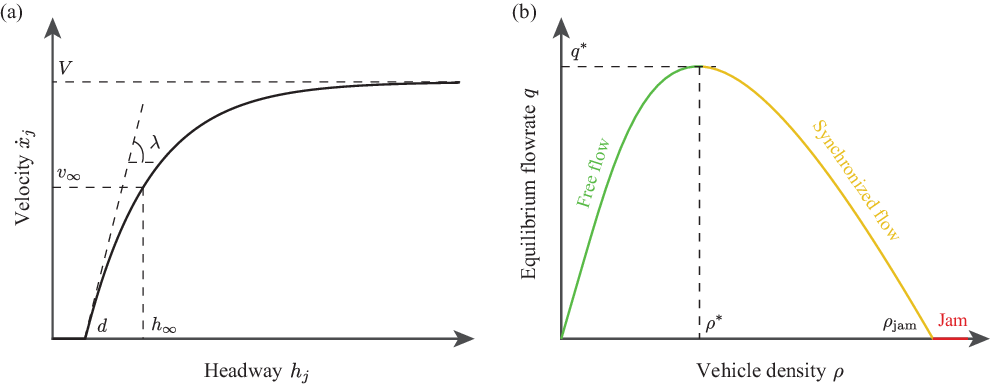}
    \caption{(a) The driving velocity $\dot{x}_j$ as a function of headway $h_j \equiv x_{j+1} - x_j$ given by \eqref{eq:velocity-eq}. $V$ is the maximal velocity, $d$ is the minimal headway for a non-zero velocity, and $\lambda$ is the rate of change of the velocity at $h_j=d$.
    In the equilibrium state where all vehicles have the same headway $h_\infty$, the equilibrium velocity of each vehicle $v_\infty$ can be found directly via \eqref{eq:velocity-eq}.
    (b) The equilibrium flow rate $q$ as a function of vehicle density $\rho$. }
    \label{fig:newell}
\end{figure}

Microscopic models describe the behavior of individual vehicles and their interactions. Microscopic models are frequently used in traffic micro simulations, which simulate the movement of individual vehicles on a road network. 
Among others, car-following models attracted a lot of attention from researchers~(See, for example, \cite{chandler1958traffic,kometani1959dynamic,newell1961nonlinear,treiterer1974hysteresis,bekey1977control,aron1988car,gipps1981behavioural,gipps1986model,ozaki1993reaction,yang1996microscopic,chen2012behavioral}). Generally speaking, car-following models postulate that the driver's behavior depends on the car (or cars) ahead of the driver. 
Many such models are second-order models that map relative positions (headway) and relative velocities to acceleration, drawing ideas from Newtonian mechanics. Such models usually involve a lot of parameters that need to be measured and calibrated. On the other hand, \citet{newell1961nonlinear} proposed a first-order model that connects the headway directly to velocity by a nonlinear mapping. Assuming identical vehicles/drivers, the velocity is given by:
\begin{equation} \label{eq:velocity-eq}
    \dot{x}_j(t) \, = \max\left( V - V \exp{\left(-\frac{\lambda}{V} [x_{j+1}(t - \Delta) - x_j(t - \Delta) - d]\right)}, 0\right)
\end{equation}
where $x_j(t)$ and $x_{j+1}(t)$ are the longitudinal positions of car $j$ and the car immediately in front of $j$ at time $t$, $x_{j+1}(t) - x_j(t)$ is the headway of car $j$, the overhead dot denotes time differentiation,  $V$ is the maximum velocity, $d$ is the minimum headway for non-zero velocity, $\lambda$ is the slope at ${x}_{j+1} - x_j(t)=d$, $\Delta \ge 0$ is the reaction time, which factors in both the driver's physiological reaction time as well as the vehicle's mechanical reaction time.
A schematic figure of \eqref{eq:velocity-eq} is shown in Figure~\ref{fig:newell}(a).

Given a traffic system with density $\rho$ (number of vehicles per unit distance), we can compute the headway of each vehicle in the equilibrium state $h_\infty$ and its corresponding equilibrium velocity $v_\infty$. The theoretical equilibrium flow rate can therefore be computed as $q = \rho v_\infty$ (number of vehicles passing a location per unit time). As shown in Figure~\ref{fig:newell}(b), this model successfully reproduces the three phases of traffic~\cite{kerner1996experimental,kerner1998experimental}. Specifically, the flow rate increases from $0$ to $q^* \equiv \max(q)$ as $\rho$ increases from $0$ to $\rho^*$, representing the free flow phase; the flow rate then decreases as $\rho$ increases to $\rho_\mathrm{jam}\equiv \arg\min_{\rho > 0} q(\rho) = 0$, representing the synchronized flow; the system enters the traffic jam phase when $\rho = \rho_\mathrm{jam}$, as the headway of each vehicle becomes smaller than the minimal headway $d$.

While Figure~\ref{fig:newell} establishes the traffic system's equilibrium state, an interesting point of research is the stability of these equilibrium states~\cite{bando1994structure,bando1995dynamical,holland1998generalised,kerner2009introduction}. That is, whether a perturbation to an equilibrium system would make the system drift away from the equilibrium state. The instability is often believed as the important contributor to the so-called phantom traffic jam: the scenario that traffic moving smoothly can suddenly spontaneously break down without obvious reasons~\cite{wilson2008mechanisms,seibold2012constructing}.
In particular, \citet{newell1961nonlinear,newell1962theories} showed that, by series expansion of the reaction time $\Delta$, the system in the continuum limit is unstable if $\lambda \Delta > 1/2$, and stable otherwise.
To better understand the source of instability in numerical simulations, \citet{kesting2008reaction} conducted a systematic study of the effects of (driver) reaction time, (vehicle) velocity adaptation time, and (simulation) update time on the stability of traffic flow. They showed that the long-wavelength string instability is mainly driven by the velocity adaptation time while short-wavelength local instabilities arise for sufficiently high reaction and update times.

Another important aspect of traffic dynamics is the modeling of lane-changing. Lane-changing is a common driving maneuver that allows drivers to navigate traffic, overtake slower vehicles, and reach their desired destinations. It introduces complexities and potential safety hazards, requiring sophisticated models to understand and predict its occurrences and impacts.
The study of lane-changing usually focuses on two different aspects: Lane Changing Decision (LCD) and Lane Changing Impact (LCI). As the name suggested, the study of LCI concentrates on analyzing the impact a lane change has on the surrounding vehicles and the overall traffic flow. It has been shown that lane changes can contribute to a reduction in roadway capacity, particularly at bottlenecks or areas of high density, and traffic oscillations~\cite{laval2006lane,laval2014parsimonious,shang2020cut,chauhan2022understanding,gao2023lane,zhang2024effects}.

On the other hand, the study of LCD focuses on the decision-making process a driver undergoes before executing a lane change. 
In his pioneer work, \citet{gipps1986model} proposed a model that contributed lane-changing decisions to factors such as gaps in the target lane for a safe lane change, the locations of permanent obstructions, the presence of transit lanes, the driver's intended turning movement, the presence of heavy vehicles, and the possibility of getting a velocity advantage. The relative importance of these factors follows a set of deterministic rules determined by the hierarchy when evaluating these factors. The hierarchy may change depending on the distance between the vehicle and its intended turning position. 
\citet{yang1996microscopic} distinguished two types of lane changes: mandatory lane changes (one which a driver must change lanes) and discretionary lane changes (one which a driver desires to change lanes). They developed and implemented a microscopic traffic simulator MITSIM that extended Gipps' deterministic model to a stochastic model that makes the model more realistic. Specifically, they introduced a probability function that approaches 1 when the vehicle is getting close to the position of a mandatory lane change. For discretionary lane changes, MITSIM considers both the driving conditions in the current lane and the adjacent lanes, and check for opportunities to improve its velocity. Several factors, such as impatient factor and velocity indifference factor, are used to determine whether the current velocity is low enough and the velocities in adjacent lanes are high enough for considering a lane change. The paper, however, does not explicitly define these factors or detail how they are implemented.
\citet{hidas2002modelling,hidas2005modelling} improved the gap acceptance criteria of \cite{gipps1986model} and classified lane changes into three categories based on the observations from video-recording footages: free, cooperative, and forced LCs. In their model, the follower's willingness to accept a certain maximum velocity decrease depends on the urgency of the lane change.
\citet{kesting2007general} proposed a lane change model that minimized overall braking induced by lane changes (MOBIL), based on the logic that the anticipated advantages and disadvantages of a potential lane change can be measured using single-lane accelerations. Furthermore, their model takes into account the potential asymmetry of lanes. For example, in countries such as Germany, the right lane is the default lane and the left lane should only be used for the purpose of overtaking, and passing in the right lane is forbidden unless traffic is congested. In these situations, MOBIL neglects the disadvantage of the immediate follower in the right lane because the left lane has priority.
Besides the Gipps-type rule-based models mentioned above, many other lane changing models exist, including the Utility-based models~\cite{ahmed1996models,ahmed1999modeling,toledo2003modeling,sun2011lane,sun2012lane}, which employ utility functions to quantify the benefits and costs associated with changing lanes or remaining in the current lane; the Cellular-automata models~\cite{nagatani1993self,rickert1996two,wagner1997realistic,nagel1998two,maerivoet2005cellular}, which model individual vehicles as particles interacting with each other, the nature of the “interactions” among these particles is determined by the way the vehicles influence each others’ movement; and the Markov-process based models~\cite{worrall1970elementary,pentland1999modeling,sheu2001stochastic,toledo2009state,singh2011estimation}, which treat lane-changing behavior as a sequence of stochastic transitions between states. 
Readers are referred to \citet{zheng2014recent} for a comprehensive review.

While significant progress has been made in developing LCD models, most of current models connect the lane-changing decisions directly to physically observable metrics such as available gaps and driving velocities. These models would work fine if all drivers are level-headed and rational. However, it is the authors' view that neglecting latent metrics such as the drivers' psychology, perception, and cognitive processes during lane changes is an over-simplification. To bridge the gap, we introduce a psychological metric, {\em frustration level} of a driver, that increases/decreases during the simulation based on the driving conditions. The frustration level, in turn, maps to the {probability} of a lane change attempt in the next second. To our knowledge, this is the first study that adds psychological components to lane-changing decision models.

The paper is organized as follows:
Section~\ref{sec:model} describes a flexible computational framework built on Newell's car-following model~\cite{newell1961nonlinear} and our proposed lane-changing mechanism. Section~\ref{sec:methods} describes the linear stability analysis for single-lane dynamics with a finite reaction time, as well as the numerical methods for the non-linear simulations. In Section~\ref{sec:results}, we show both analytically and numerically that there is a critical reaction time that depends on the vehicle density, beyond which the vehicles will collide. Additionally, we show that our lane-changing model is capable of reproducing real-life phenomena such as load-balancing between lanes. We also examine the effects of being an ``aggressive'' driver on a busy road. Finally, we discuss the limitations and future directions of the work in Section~\ref{sec:discussion}.

\section{Model}
\label{sec:model}
\begin{figure}
    \centering
    \includegraphics[]{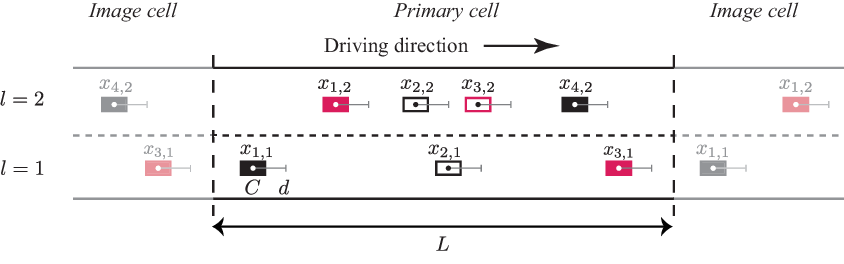}
    \caption{Schematic figure of our car-following model. Each rectangle denotes a vehicle whose position $x_{j, l}$ is represented by the dot at the center of the rectangle. The width of each rectangle represents the size of the vehicle $C$, and the length of the gray whisker denotes the minimal headway $d$. Red rectangles denote the vehicles whose drivers are accumulating frustration because the headway in the adjacent lane is bigger than the headway in its own lane. Filled rectangles denote the vehicles that have enough gap in the adjacent lane to complete a lane change (no car in [$x_{j,l} - d, x_{j,l} + d$]), while open rectangles do not. The simulation is done in a primary cell of length $L$. All positions are measured from the left end of the primary cell. Vehicles in the image cells are shown in lighter colors, as a reminder that they are merely images of the vehicles of the primary cell. }
    \label{fig:model}
\end{figure}

Consider a total of $N$ vehicles in multiple lanes as shown in Figure \ref{fig:model}. $x_{j,l}$ denotes the position of the $j$-th vehicle in the $l$-th lane. 
The headway of each vehicle is measured from the center of itself to the center of the vehicle in front of it $h_{j,l} = x_{j+1,l} - x_{j,l}$, and this headway is used to compute the driving velocity of this vehicle via \eqref{eq:velocity-eq}.
We adopt a periodic boundary condition such that there is no ``leading'' vehicle as each vehicle has one in front of them. For example, vehicle 3 is following vehicle 1 in lane 1 and the headway is $h_{3,1} = x_{1,1} + L - x_{3,1}$, where $L$ is the size of the primary simulation cell.
The size of each vehicle is denoted by the width of the rectangle $C$ to determine whether there is a collision in the system. Cars $j$ and $j+1$ are considered in collision if $h_{j,l} = x_{j+1,l} - x_{j,l} \le C$.
It is important to make sure that the minimal headway $d$ in \eqref{eq:velocity-eq} is greater than the vehicle size so that the vehicle tries to stop before running into the vehicle in front of it.

We focus on discretionary lane changes in this work rather than mandatory lane changes. Our underlying hypothesis is that drivers want to achieve the fastest possible velocity according to \eqref{eq:velocity-eq} by making lane changes to position themselves in the ``best'' lane. To facilitate the decision, we introduce the {\em frustration level} $\phi(t)$ of a driver to compute their probability of a lane change attempt at time $t$. Given two parameters $r\ge 0$ and $p\ge 0$, the frustration level $\phi(t)$ changes according to the following four rules in our model:
\begin{enumerate}
    \item ${d\phi}/{dt} = r$ if the driver's headway in the current lane is smaller than the perceived headway in an adjacent lane.
    \item ${d\phi}/{dt} = -r$ if the driver's headway in the current lane is larger than the perceived headway in an adjacent lane. 
    \item $\phi$ will increase by $p$ every time the driver is passed by one vehicle in an adjacent lane.
    \item $\phi$ will reset to 0 after each lane change.
\end{enumerate}
Simply put, the driver's frustration level will increase if they believe they are currently in the ``worse'' lane compared to the adjacent lane, and decrease if they believe they are in the ``better'' lane. 
In terms of Figure~\ref{fig:model}, the drivers in the vehicles labeled in red color will raise their frustration levels based on rule \#1, and the other vehicles will reduce their frustration levels based on rule 2.
Note that we keep $\phi(t)\ge 0$ in all cases regardless of rule \#2.

We propose a stochastic model with an underlying function $P(\phi)$ that maps the frustration level to the probability of a lane change attempt in the next second
\begin{equation}
    P(\phi) = \frac{2}{\pi} \arctan(\phi).
\end{equation}
We note that there is no best way to map the frustration level to the probability of a lane change attempt without the support of empirical data. 
However, this simple function satisfies the following three important properties to serve as the probability of a lane change attempt:
\begin{enumerate}
    \item $0 \le P(\phi) \leq 1$ for all $\phi \in \mathbb{R}^{\ge 0}$,
    \item $P(\phi)$ is monotonically increasing,
    \item $\lim_{\phi \rightarrow \infty}{P(\phi)} = 1$.
\end{enumerate}
We stress that $P(\phi)$ is not the probability of a lane change but a lane change {\em attempt}. The attempt will be successful if there is no car in the adjacent lane in the interval $[x_j - d, x_j + d]$, and unsuccessful otherwise. In the event of an unsuccessful lane change attempt, the vehicle will stay in its original lane. 
For example, in Figure~\ref{fig:model}, the lane change attempt of vehicle 1 in lane 2 will be successful, while the attempt of vehicle 2 in lane 2 will be unsuccessful. 
It is also important to note that a lane change attempt can happen even when the driver's frustration level is decreasing as long as $P(\phi) > 0$ (for example, vehicle 1 in lane 1). As a result, not all lane changes will increase the lane-changer's headway.

\section{Methods}
\label{sec:methods}
\subsection{Linear Stability Analysis}
\label{sec:LSA}
The linear stability of $N$ identical vehicles in a {\em single-lane} can be examined analytically. 
Let $y_j(t) = x_j(t) - x_N(t)$ be the relative position of the $j$-th vehicle with respect to the $N$-th vehicle. Note that we have omitted the lane index here.
Ignoring the trivial case of $j=N$ in which $y_N(t) \equiv 0$, we can get a vector ${\mathbf y}(t) = \{y_j(t)\}$ for $j = 1, 2, \cdots, N-1$. 
Differentiating $\mathbf{y}$ with respect to time, we get a vector representing the relative velocity of the vehicles $\dot{{\mathbf y}}(t) = \{\dot{y}_j(t)\} = \{\dot{x}_j(t) - \dot{x}_N(t)\}$, where $\dot{x}_j(t)$ and $\dot{x}_N(t)$ are given by equation~\eqref{eq:velocity-eq}.
Specifically, 
\begin{equation}\label{eq:relative-vel}
\begin{split}
    \dot{y}_j(t) &= \dot{x}_j(t) - \dot{x}_N(t)\\
    &= -V\exp\left(-\frac{\lambda}{V}[x_{j+1}(t-\Delta) - x_j(t-\Delta) - d]\right) + V\exp\left(-\frac{\lambda}{V}[x_1(t-\Delta)+L - x_N(t-\Delta) - d]\right)\\
    &= -V\exp\left(-\frac{\lambda}{V}[y_{j+1}(t-\Delta) - y_j(t-\Delta) - d]\right) + V\exp\left(-\frac{\lambda}{V}[y_1(t-\Delta)+L - d]\right),
\end{split}
\end{equation}
where the track length $L$ is added in the second term to account for the periodic boundary condition. 
Linearizing equation~\eqref{eq:relative-vel} about the equilibrium state in which all vehicles have the same headway $L/N$, we have 
\begin{equation}\label{eq:linearized}
    \dot{{\mathbf y}}(t) = \mathbf{J}\,{\mathbf y}(t-\Delta),
\end{equation}
where $\mathbf{J} = \{J_{jk}\} = \{ \partial \dot{y}_j/\partial y_k\}$ is the Jacobian matrix. 
Due to the nature of the car-following model, the Jacobian matrix is a sparse matrix that can be easily expressed as $\mathbf{J} = c(\mathbf{I}+\mathbf{A}+\mathbf{B})$, where $c=-\lambda \exp{\left(-\frac{\lambda}{V} (\frac{L}{N} - d)\right)}$ is a (negative) scaling constant, $\mathbf{I}$ is the identity matrix, 
 the matrices $\mathbf{A}$ and $\mathbf{B}$ can be described as 
\begin{equation}\label{eq:AB}
A_{jk} = 
\begin{cases} 
-1 & \text{if } k = j + 1 \\ 
0 & \text{otherwise}
\end{cases},
\qquad \text{ and }\qquad
B_{jk} = 
\begin{cases} 
1 & \text{if } k = 1 \\ 
0 & \text{otherwise}
\end{cases}.
\end{equation}

The characteristic equation of the linearized delay differential equation \eqref{eq:linearized} is, by \citet[Section 4.3]{smith2011introduction}:
\begin{equation}\label{eq:characteristic}
    \det(\tilde\lambda \mathbf{I} - \mathbf{J} \exp(-\tilde\lambda \Delta)) = 0,
\end{equation}
and the real parts of the characteristic roots of this equation $\mathrm{Re}(\tilde{\lambda})$ describe the stability of the dynamical system. 
Naively solving this equation is extremely computationally expensive for a moderate number $N$ and nearly impossible for large $N$. 
However, we can pursue a much more efficient algorithm based on the fact that $\mathbf{J}$ is diagonalizable (proof at the end of the sub-section). That is, we can find an invertible matrix $\mathbf{R}$ such that 
$\mathbf{J} = \mathbf{RDR}^{-1}$, 
 where $\mathbf{D}$ is the diagonal matrix formed by the eigenvalues of $\mathbf{J}$. Following this, we can simplify the left-hand side of equation~\eqref{eq:characteristic} as such:
\begin{equation}
     \det(\tilde\lambda \mathbf{I} - \mathbf{J} \exp(-\tilde\lambda \Delta)) 
    =\det(\mathbf{R}\tilde\lambda \mathbf{I R}^{-1} - \mathbf{RDR}^{-1} \exp(-\tilde\lambda \Delta)) 
    = \det(\mathbf{R}(\tilde\lambda \mathbf{I} - \mathbf{D} \exp(-\tilde\lambda \Delta))\mathbf{R}^{-1}) 
    =\det(\tilde\lambda \mathbf{I} - \mathbf{D} \exp(-\tilde\lambda \Delta)).
\end{equation}
Since $\tilde\lambda \mathbf{I} - \mathbf{D} \exp(-\tilde\lambda \Delta)$ is a diagonal matrix, the determinant is simply the product of all diagonal elements. 
The characteristic equation~\eqref{eq:characteristic} can therefore be rewritten as
\begin{equation}\label{eq:characteristic2}
    \prod_{j}\left[\tilde\lambda - \tilde d_j \exp(-\tilde\lambda \Delta)\right] = 0,
\end{equation}
where $\tilde d_j$ is the $j$-th eigenvalue of the Jacobian $J$. We use a Python library {SymPy} to solve the characteristic roots $\tilde\lambda$.

\vskip 5pt
\noindent
\underline{Proposition}: The Jacobian matrix in equation~\eqref{eq:linearized} is diagonalizable. 
\vskip 5pt
\noindent
\underline{Proof}: Consider the $(N-1)\times (N-1)$ square matrix $\mathbf{M} = \mathbf{J}/c = \mathbf{I} + \mathbf{A} + \mathbf{B}$, where $\mathbf{I}$ is the identity matrix and $\mathbf{A}$ and $\mathbf{B}$ are defined in \eqref{eq:AB}. We will show that $\mathbf{M}$ has $N-1$ distinct eigenvalues. The characteristic polynomial of $\mathbf{M}$, $\mathcal{P}(\tilde\lambda) = \det(\mathbf{M} - \tilde\lambda \mathbf{I})$, can be obtained by Laplace expansion along the first row of the matrix $\mathbf{M}-\tilde\lambda \mathbf{I}$: $\mathcal{P}(\tilde\lambda) = (2 - \tilde\lambda)\mathcal{P}_{(1, 1)}(\tilde\lambda) + (-1)^2 \mathcal{P}_{(1, 2)}(\tilde\lambda)$, where $\mathcal{P}_{(j, k)}(\tilde\lambda)$ is the ($j$, $k$) minor of matrix $\mathbf{M} - \tilde\lambda \mathbf{I}$.  
The (1, 1) minor can be very easily computed as the corresponding matrix is diagonal: $\mathcal{P}_{(1, 1)} = (1 - \tilde\lambda)^{N - 2}$. Continuously applying Laplace expansion along the first row of the corresponding matrix of the (1, 2) minor yields $\mathcal{P}_{(1, 2)} = (1 - \tilde\lambda)^{N - 3} + (1 - \tilde\lambda)^{N - 4} + \cdots + 1$. Therefore, 
\begin{equation*}
    \begin{split}
 \mathcal{P}(\tilde\lambda) &=(2 - \tilde\lambda)\mathcal{P}_{(1, 1)}(\tilde\lambda) + \mathcal{P}_{(1, 2)}(\tilde\lambda)\\
 &=(1 + (1 - \tilde\lambda))(1 - \tilde\lambda)^{N - 2} + (1 - \tilde\lambda)^{N - 3} + (1 - \tilde\lambda)^{N - 4} + \cdots + 1 \\
 &= (1 - \tilde\lambda)^{N-1} + (1 - \tilde\lambda)^{N - 2} + (1 - \tilde\lambda)^{N - 3} + (1 - \tilde\lambda)^{N- 4} + \cdots + 1\\
 &=\sum_{k=0}^{N-1} (1-\tilde\lambda)^k.
    \end{split}
\end{equation*}
The solutions for $\mathcal{P}(\tilde\lambda) = 0$ are given by $\tilde\lambda = 1 - \exp(2 \pi \mathrm{i}k/N)$, where $\mathrm{i} = \sqrt{-1}$ and $k = 1, 2, 3, \cdots, N-1$. Therefore, the matrix $\mathbf{M}$ has $N-1$ distinct eigenvalues and is diagonalizable, so is $\mathbf{J} = c\mathbf{M}$.

\subsection{Numerical Method}
\label{sec:numerical}
While the linear stability analysis can provide us invaluable insights into the dynamical system of single-lane, we need a numerical method to study the case of multi-lane systems. Additionally, a numerical method allows us to probe into each simulation and study the spatiotemporal feature embedded in the dynamical system.

\begin{figure}
    \centering
\includegraphics[scale = 0.95]{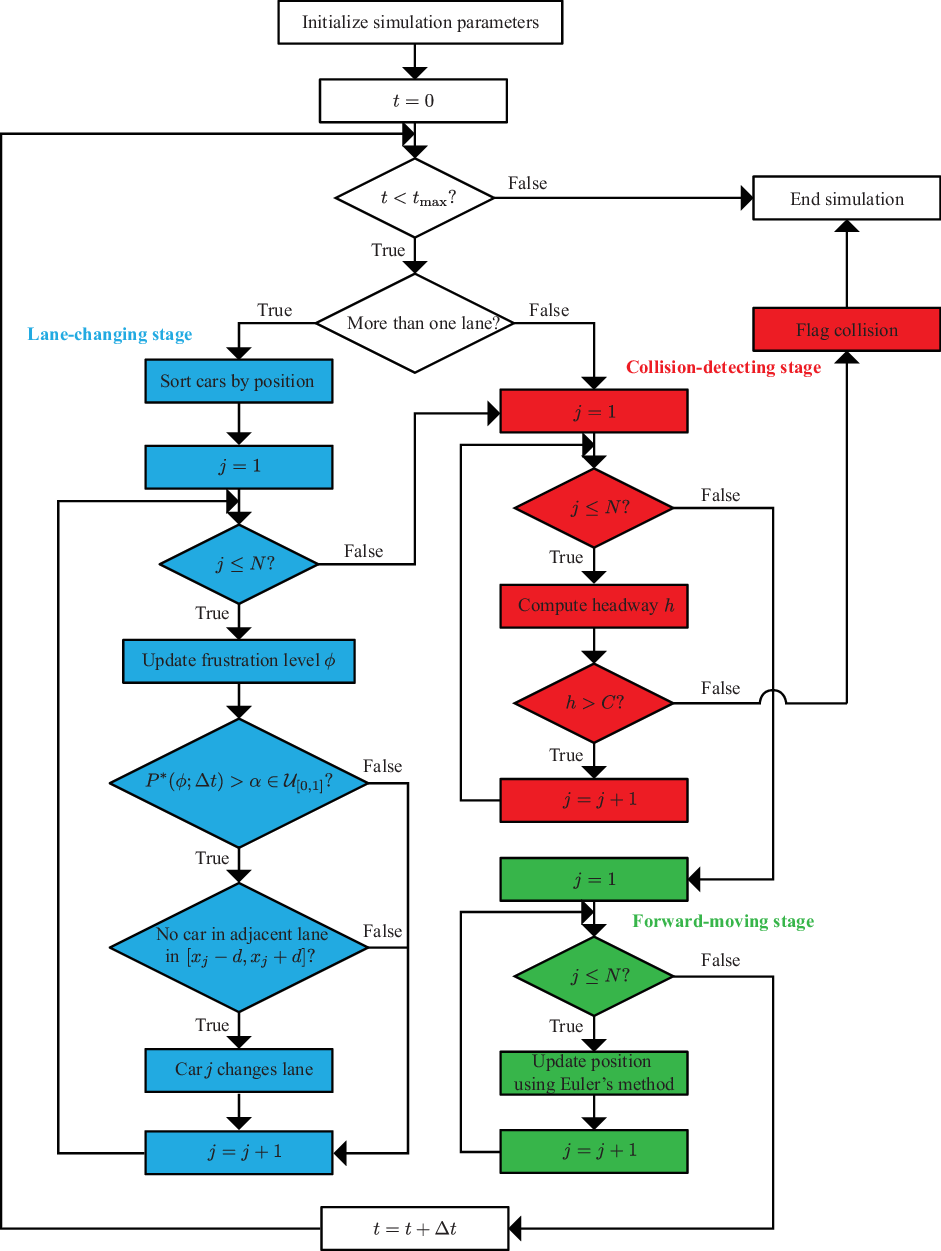}
\caption{Flowchart of the numerical method. The steps are grouped into Lane-changing, Collision-detecting, and Forward-moving stages, denoted by blue, red, and green colors respectively.}
\label{fig:flowchart}
\end{figure}

The algorithm is summarized as a flowchart in Figure~\ref{fig:flowchart}. 
At each time step, the simulation goes through a lane-changing stage, a collision-detecting stage, and a forward-moving stage. In the lane-changing stage, the algorithm loops through each vehicle $j$, updates the frustration level based on the four rules listed in Section~\ref{sec:model} and decides whether to make a lane-change attempt and subsequently whether such an attempt is successful. In the collision-detecting stage, the algorithm loops over each vehicle and detects whether there is a collision. Finally, in the forward-moving stage, each vehicle's position is updated using Euler's method with a fixed time step $\Delta t$, where the velocities are given by \eqref{eq:velocity-eq}. 
Note that the headway in the lane-changing stage and forward-moving stage are computed using the {\em perceived} headway (taken into account the reaction time $\Delta$), while the collision-detecting stage uses the {\em true} headway.

A short explanation of the lane-changing stage is provided. Particularly, in addition to its own position and its headway, each vehicle also keeps track of the position of the vehicle immediately in front of itself in the adjacent lane to determine the headway in the adjacent lane, which is used to determine whether the driver's frustration level will increase or decrease;
furthermore, the position of the vehicle immediately behind the vehicle in the adjacent lane is also used, in addition to the vehicle in front, to determine whether a lane change would be safe.
The frustration level $\phi$ is updated at each time step for each driver. To decouple the choice of time step and the probability of the lane-change attempt, we rescale the lane change probability as
\begin{equation}\label{eq:rescaled-prob}
    P^*(\phi;\Delta t) = 1 - [1-P(\phi)]^{\Delta t}.
\end{equation}
A random number is drawn from a uniform distribution $\mathcal{U}_{[0, 1]}$ at each time step for each vehicle and compared with $P^*(\phi; \Delta t)$ to determine whether the vehicle would attempt a lane change at this time step. 
The vehicle's longitudinal position is not affected by the lane change.
It is easy to verify that the rescaled probability per time step recovers the probability of a lane change attempt per unit time by considering the probability of lane change attempt in the consecutive $1/\Delta t$ time steps: $1 - [1 - P^*(\phi;\Delta t)]^{1/\Delta t} = P(\phi)$.

The parameters used in the remaining of our paper are listed in Table~\ref{tab:parameters} unless otherwise stated. The maximal equilibrium flow rate given by these parameters is $q^* = 2065$ vehicles per hour, corresponding to the critical vehicle density $\rho^* = 34$ vehicles per kilometer. Additionally, $\rho_\mathrm{jam} = 134$ vehicles per kilometer. These values are in good agreement with prior experimental results \cite{kerner1996experimental,kerner1998experimental}.
\begin{table}[H]
\centering
\begin{tabular}{|
>{\columncolor[HTML]{C0C0C0}}c |c|c|c|c|c|}
\hline
{\color[HTML]{000000} \textbf{Parameter}} & Velocity change rate                           & Maximal velocity          & Minimal headway & Vehicle size   & Track length     \\ 
\hline
{\color[HTML]{000000} \textbf{Symbol}} & $\lambda$                           & $V$          & $d$ & $C$   & $L$     \\ 
\hline
\textbf{Value}                            & \cellcolor[HTML]{FFFFFF}1.0~/s & 40~m/s & 7.5~m           & 5~m    & 1000~m               \\ \hline
\end{tabular}
\caption{Physical parameters of the simulations. }\label{tab:parameters}
\end{table}

\subsection{Metrics}
\label{sec:metrics}

In this subsection, we list all numerical metrics that we use in the subsequent Results section.

To numerically study the stability of the single-lane system, we look at the {\em cyclic growth rate of velocity difference}.  In particular, define a function $f: \mathbb{N} \rightarrow \mathbb{R}$ that represents the amplitude of the velocity oscillation as follows:
\[
    f(n) = v_\infty - v_{\mathrm{min}, n},
\]
where $v_\infty$ is the equilibrium velocity computed by \eqref{eq:velocity-eq} assuming equal headway for each vehicle, and $v_{\mathrm{min, n}}$ is the minimum velocity of the $n$th period for one vehicle.  
Fitting an exponential curve $ae^{kn}$ to $f(n)$, the value of $k$ is the cyclic growth rate of velocity difference. The system is stable if $k<0$, and unstable otherwise.

To measure the capacity of the lane(s), we compute the flow rate $Q$ at time $t$ and position $x$ by
\begin{equation}\label{eq:flowrate}
      Q(t, x) = \frac{1}{\delta } \sum_{j = 1}^N [H(x_j(t) - x) - H(x_j(t-\delta) - x)],
\end{equation}
where $H(x):=\left\{\begin{array}{ll}1, & x\ge 0\\ 0, & x < 0\end{array}\right.$ is the Heaviside step function, and $\delta > 0$ is the (numerical) averaging time window. In practice, the time window $\delta$ is chosen large enough such that there is at least several vehicles passing the position $x$ within the window, and small enough such that no vehicle would drive the entire track.
Intuitively, the flow rate $Q(x, t)$ measures the number of vehicles passing a particular longitudinal position $x$ at time $t$ per time unit. 
In the case of multi-lane simulations, the flow rate is summed over all lanes.

Additionally, in the case of multi-lane simulations, we measure the absolute difference in the number of vehicles in each lane $\Delta N$ at any given time to understand the load balance of the lanes. We also keep track of the total number of lane changes for each vehicle $\Delta l$ during the simulation.

\section{Results}
\label{sec:results}
\subsection{Single-lane dynamics}
\begin{figure}
    \centering
    \includegraphics[width=\textwidth]{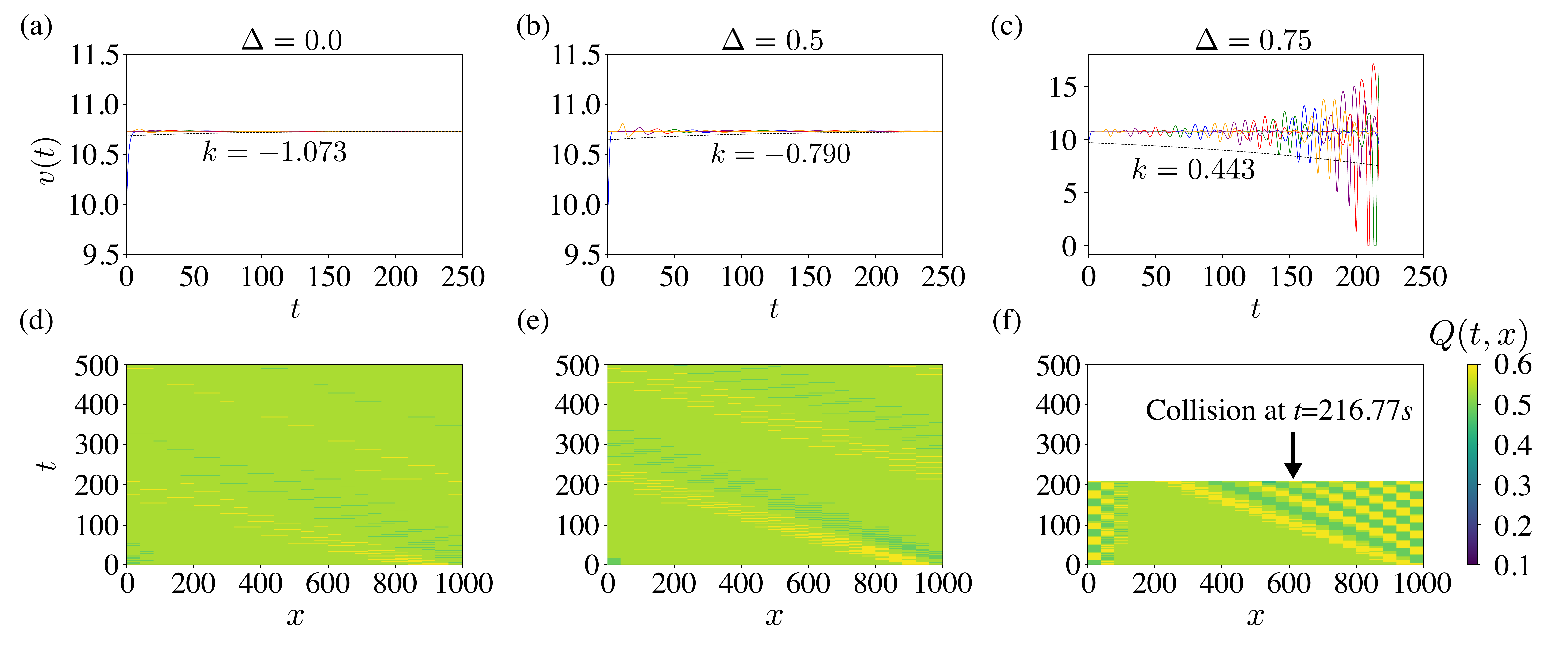}
    \caption{Dynamics of 50 vehicles in a single lane with various reaction times $\Delta = 0, 0.5, 0.75$s. (a-c) Velocity as a function of time for representative vehicles with different reaction times. Curves of blue, orange, purple, red, and green represent the results of vehicles 1, 11, 21, 31, and 41 respectively. Dashed curves denote the decay/growth of velocity difference.
    (d-f) The flow rate at different position $x$ and time $t$ for various reaction time $\Delta$. The time step and the averaging time window for all simulations are $\Delta t = 0.01s, \delta = 18.63s$.}
    \label{fig:wave-visualizations}
\end{figure}

We start by asking the question of how the reaction time $\Delta$ affects the stability of single-lane dynamics by examining how the system responds to a perturbation of initial positions with different reaction times. 
Specifically, we place $N= 50$ vehicles uniformly along the track such that $x_j(0) = (j-1)L/N$ for $j \in \{1, \cdots, N\}$. This is the system's equilibrium state and all vehicles will move at the same constant velocity $v_\infty$ regardless of the reaction time. 
The equilibrium velocity $v_\infty$ can be easily computed using equation~\eqref{eq:velocity-eq} by setting $x_{j+1}(t) - x_j(t) = L/N$. 
We introduce a perturbation to the system by displacing vehicle $j=1$ by 1~m downstream. This perturbation decreases the initial velocity of vehicle $j=1$ and increases the initial velocity of vehicle $j=N$. This initial perturbation also introduces velocity oscillations for each vehicle as they evolve according to \eqref{eq:velocity-eq}. The velocities of representative vehicles with $\Delta$ = 0, 0.5~s, and 0.75~s are shown in Figure~\ref{fig:wave-visualizations}(a-c).
In the cases of $\Delta = 0$ and 0.5~s, 
the amplitude of the oscillations decreases over time as all vehicles converge to the velocity of the equilibrium state $v(t) \equiv \dot{x}(t) \rightarrow v_\infty$ as $t\rightarrow \infty$. The  convergence rate is lower in the case of $\Delta = 0.5$~s, as $k$ increases from $-1.073$ to $-0.790$.
In the case of $\Delta = 0.75$~s, however, the differences in velocities increase as time progresses, resulted in a positive $k = 0.443$. 
A crash is identified at $t\approx 217$s and the simulation is ended.

To further understand the spatiotemporal features of the traffic flow, we plot the flow rate $Q(t,x)$ in Figure~\ref{fig:wave-visualizations}(d-f). 
The equilibrium flow rate $q$ of the system is $q = 0.54$ vehicles/second ($1944$ vehicles/hour).
Generally speaking, a flow rate lower than $q$ indicates a congestion area, indicated by blue color in the panels. In all cases, the congestion areas propagate ``backward'' (goes upstream) as time progresses, consistent with many previous observations and models~(See, for example, \citet{kerner1999physics,bando1994structure,chen2014periodicity}). In the zero and small reaction time cases (Fig~\ref{fig:wave-visualizations}(d,e)), the congestion is dissipated over time, and the rate of dissipation with $\Delta = 0.5$~s is slower, following the same trend of the rate of convergence of velocities seen in Figures~\ref{fig:wave-visualizations}(a-b). 
On the other hand, the congestion area does not show the self-dissipating capability in the case of $\Delta = 0.75$~s. Instead, the initial congestion area produces multiple congestion area feeding into each other, indicated the alternating light blue-yellow patterns at given time.

\begin{figure}
    \centering
    \includegraphics[width = 0.9\textwidth]{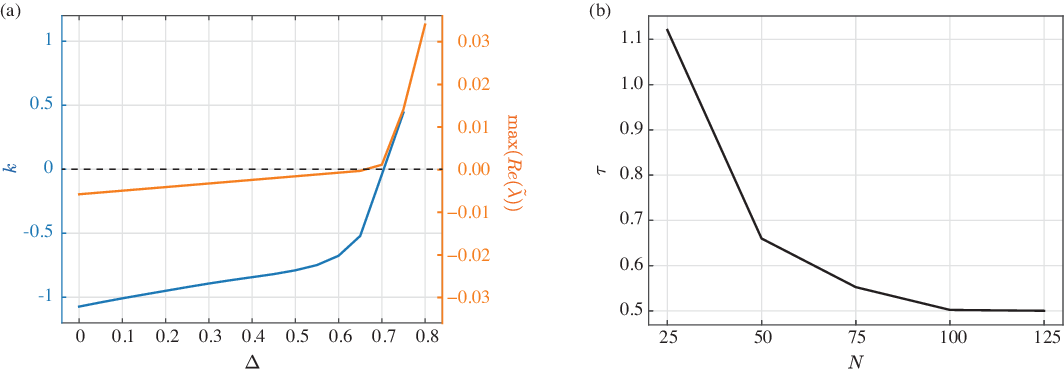}
    \caption{Stability of a single-lane traffic system. (a) The cyclic growth rate of velocity difference $k$ of the non-linear simulations and the maximal real part of the characteristic roots of the linearized system as functions of reaction time $\Delta$. (b) The critical reaction time $\tau$ as a function of the number of vehicle in the system $N$.}
    \label{fig:stability}
\end{figure}

These results lead us to believe that there exists a critical value of reaction time $\tau$ such that the system is stable if $\Delta < \tau$, and unstable otherwise.
To this end, we plot the cyclic velocity difference growth rate $k$ as well as the maximal real part of the characteristic roots of equation~\eqref{eq:characteristic} $\max_j (\mathrm{Re}(\tilde\lambda_j))$ against the reaction time $0\le\Delta\le 0.8$~s at a fixed interval 0.05~s in Figure~\ref{fig:stability}. 
We note that the simulation detects a collision within one cycle for $\Delta  = 0.8$s, prohibiting us from obtaining a $k$ value.
Our results show that both the linear stability analysis and the non-linear simulations reveal that a longer reaction time is more likely to result in an unstable system. Additionally, both models predict a critical reaction time $\tau\approx 0.7$~s.

Linear stability analysis is also conducted on cases of different $N$, shown in Figure~\ref{fig:stability}(b).
The critical reaction time $\tau$ decreases with the number of vehicles $N$, consistent with our intuition as higher vehicle density reduces the space between adjacent vehicles, leaving shorter time to ``react''.
We note that $\tau$ approaches 0.5~s as the system approaches its theoretical limit $N_\mathrm{max} = L/d = 1000/7.5 \approx 133$. This case can be considered as the ``continuum limit'', and our result recovers the classic result shown in \citet{newell1961nonlinear}.

\subsection{Multi-lane dynamics}

We next showcase the effects of our lane-changing mechanism. Specifically, we demonstrate the results of two experiments that displays a natural ``load-balancing'' behavior and the how adding an ``aggressive driver'' can disrupt the system. In both experiments, we set the reaction time $\Delta = 0$ to focus on the proposed lane-changing mechanism.

\begin{figure}
    \centering
    \includegraphics[width=0.6\linewidth]{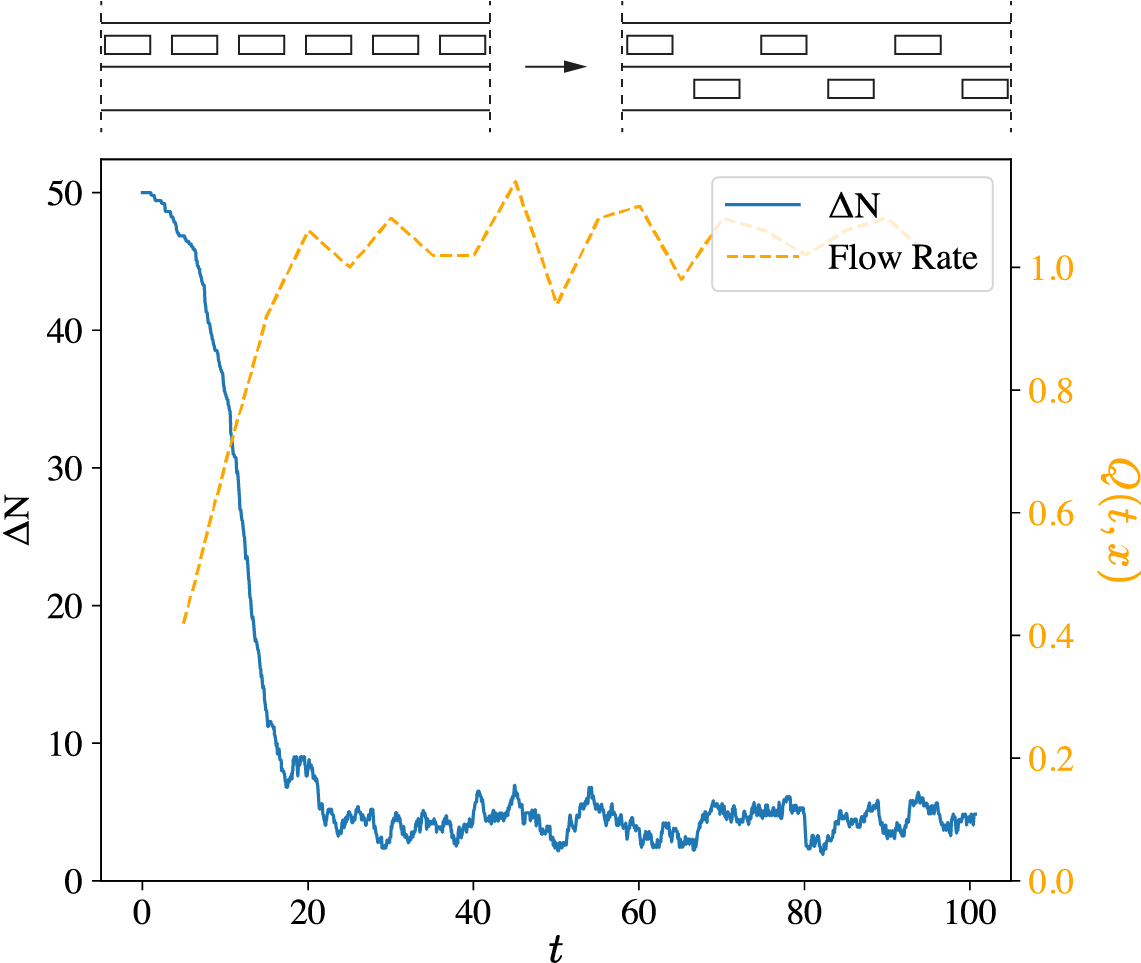}
    \caption{
    The absolute difference between the numbers of vehicles in each lane $\Delta N$ and the flow rate $Q(t, L/2)$ as functions of time.
    The curves are the mean value of 10 Monte Carlo simulations.
    The averaging time window for flow rate calculation is $\delta = 5$s.  All vehicles are initially placed in one lane with equal spaces.
    $N = 50$, $\Delta t = 0.05$s, $r = 0.1$, $p = 0.2$.
    }
    \label{fig:load-balance}

\end{figure}

To demonstrate the results of load-balancing, we initialize the simulation such as all vehicles are in one lane with equal spaces and the other lane is empty. 
Absence of a lane-changing mechanism, all vehicles will keep moving in the original lane as the equal headway establishes an equilibrium state. With our lane-changing mechanism in place, all driver' frustration levels start to increase at rate $r$ as the headway in the adjacent lane ($\infty$) is greater than its current headway. 
Note that not all vehicles make the lane-changing attempt at the same time because of our stochastic model, even though the frustration levels increase at the same rate at the beginning. Specifically, suppose all vehicles have the same frustration level $\phi = 0.5$ at some time $t$ in the simulation, the probability of making a lane-changing attempt for each vehicle in the next second is then $P(\phi) = \frac{2}{\pi}\arctan(0.5) \approx 30\%$. Looking at the system as a whole, this means that about 30\% of all vehicles will make a lane-changing attempt within the next second. 
To understand the time evolution of the system, we trace the absolute difference in the number of vehicles between the two lanes $\Delta N$ as well as the flow rate at the middle of the primary cell $Q(t, L/2)$ over the entire simulation. 
To reduce the random effect in the results, we conduct 10 Monte Carlo simulations and plot the average $\Delta N$ of these 10 runs for the first 100 seconds in Figure~\ref{fig:load-balance}.
Evidently, $\Delta N$ gradually reduces from 50 to about 5 within about 30 seconds and keeps oscillating around 5. 
The gradual decrease of $\Delta N$ recovers what we observe in the real-world, as an empty lane will soon be occupied by as many vehicles as in the adjacent lane to make both lanes ``symmetric''.
We note that perfectly balanced lanes $\Delta N = 0$ with equal headway is an equilibrium state. Our simulations, however, do not seem to reach that equilibrium state even after a time much longer than $100$s, likely due to the choice of our parameters $r$ and $p$.
Smaller values of $r$ and $p$ would potentially reduce both the rate of the decrease of $\Delta N$ as well as $\Delta N$ at large $t$, as the system would get more time to ``relax''. 
Additionally, the average flow rate measured at the middle of the simulation cell is also plotted in Figure~\ref{fig:load-balance} and shows an opposite trend as $\Delta N$. Specifically, the flow rate increases with time, and stagnates at around 1.0 vehicle/second (or 3600 vehicles per hour) after about 30 seconds.
Note that this is the flow rate for {\em both} lanes, and each lane has about 1800 vehicles per hour on average.
To put the result in perspective, the theoretical equilibrium flow rate with 25 vehicles per kilometer is computed as $q = 1668$ vehicles per hour, and the maximal equilibrium flow rate $q^* = 2065$.
Our result shows that the proposed lane-changing mechanism automatically approaches a globally optimal state in which the number of vehicles in each lane is roughly the same and maximizes the flow rate of the system.

Next, we test our model's capacity by introducing inhomogeneity into the system.
Specifically, we increase one driver's $\lambda$ while keeping all other parameters the same as in Table~\ref{tab:parameters}. 
Compared to other drivers, increasing $\lambda$ has two direct effects:
\begin{enumerate}
\item The driver's velocity is more sensitive to the headway when $h\approx d$.
\item The driver would keep a smaller headway to drive at the same velocity.
\end{enumerate}
As a result of our lane-changing mechanism, effect 2 makes the driver find themselves in a place where the frustration level will be raised. This leads to a third effect of increasing $\lambda$:
\begin{enumerate}
\setcounter{enumi}{2}
\item The driver tends to attempt lane-changes more frequently than the other drivers.
\end{enumerate}
These characteristics make the driver with a higher $\lambda$ a suitable case for a typical ``aggressive driver''~\cite{james2009road}.
In the following, we refer to the driver with $\lambda = 2$ as the {aggressive driver} and all other drivers as {control drivers} who have $\lambda = 1$.

As the initial condition, we place 25 vehicles in each lane with equal headway. Additionally, the vehicles in lane 1 are shifted downstream by $L/50$ so that the vehicles in both lanes are ``staggered'' similar to the configuration shown in the upper-right panel of Figure~\ref{fig:load-balance}. This is an equilibrium state if all drivers are identical as they all have the same headway, and the headway is bigger than the headway in the adjacent lane so that they do not raise their frustration levels. Each vehicle will drive at the same velocity and execute zero lane changes.

\begin{figure}
    \centering
    \includegraphics[]{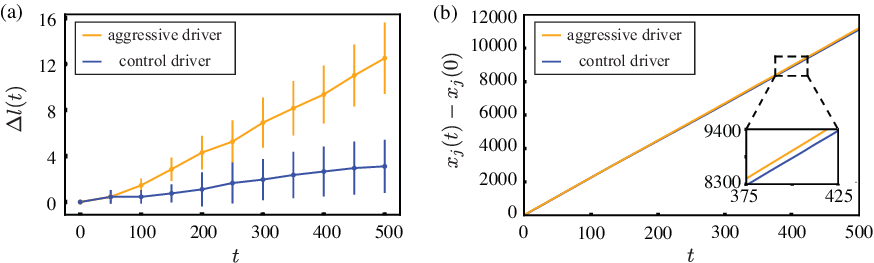}

    \caption{(a) Number of lane-changes $\Delta l$ of the aggressive driver ($\lambda = 2$) and a control driver ($\lambda = 1$) over time. Curve and the error bars shows the mean and the standard deviation of 20 Monte Carlo simulations.
    (b) Driving distances of the aggressive driver and a control driver over time. 
    $N = 50$, $\Delta t = 0.05$s, $r = 0.1$, $p = 0.1$.}
    \label{fig:aggressive}
\end{figure}

To examine how an aggressive driver affects the traffic system, we make the driver of one vehicle in lane 1  the aggressive driver. Figure~\ref{fig:aggressive}(a)  compares their number of lane changes over time against that of a control driver initially in the next lane.
The aggressive driver executes about 4 times as many lane changes as the control driver on average over 20 Monte Carlo simulations (12 times in 500 seconds vs 3 times in 500 seconds). 
The result also shows that the introduction of the aggressive driver increases the number of lane changes of the control driver as the system is being pushed away from its equilibrium state.
It is worthwhile to point out that in our model, the aggressive driver themselves does not immediately raise their frustration level as the initial positions are staggered. However, the aggressive driver keeps a smaller headway compared to all other drivers, effectively giving other drivers in lane 1 a bigger average headway than the ones in lane 2 and making all drivers in lane 1 drive faster. The difference in velocity would eventually make the aggressive driver pass one vehicle in lane 2 and start to raise their frustration level and induce a jump of the frustration level for the driver being passed.
On the other hand, our simulation also shows that there is little advantage of being an aggressive driver in terms of driving distance.
Specifically, Figure~\ref{fig:aggressive}(b) shows that the aggressive driver would drive at a velocity less than $2\%$ faster than the control driver.

\section{Discussion}
\label{sec:discussion}
In this work, we proposed a computational framework based on a car-following model capable of analyzing traffic stability and lane changes.
Specifically, we found a critical reaction time $\tau$ that depends on the vehicle density in single-lane traffic. The critical reaction time determines whether the system is stable to an initial perturbation to the equilibrium system. 
We introduced a stochastic lane-changing mechanism that bridges a driver's psychological measure and the probability of a lane-change attempt. 

Previous studies established that traffic oscillations strongly correlate with collisions~\cite{zheng2010impact}.
In our stability analysis, we examine how the traffic system responds to small oscillations initially induced by a small perturbation to the equilibrium state. 
We find that the initial oscillation decays over time if the reaction time of the drivers is low. On the other hand, the initial oscillation grows into multiple oscillations, and eventually leads to a collision if the reaction time of the drivers is high. The critical reaction time that distinguishes stable and unstable systems is studied both analytically and numerically and the results are in good agreement.
In particular for our analytical study, the problem is formulated as a system of delay-differential equations. We developed an efficient algorithm to find the characteristic roots of the linearized system to determine the stability of the system.
The critical reaction time reduces with the vehicle density $\rho$, and converge to the classical result of \citet{newell1961nonlinear} in the continuum limit.

We proposed a psychological measure of the driver in our lane-changing model, {\em frustration level}, that maps to the {\em probability} of a lane change attempt. Four explicit rules govern how the frustration level changes during the simulation. Generally speaking, the frustration level increases if the driver feels that they are in the ``worse'' lane than the adjacent lane, and decreases otherwise. We demonstrated that the proposed mechanism can reproduce interesting real-life phenomena like load-balancing of multiple lanes. We also applied our model to study the effects of being an aggressive driver and showed that while being an aggressive driver yields a 2\% time benefit, it costs four times as many lane changes as a control driver, posing significant risks to road safety~\cite{zheng2013effects}.


There are many extensions that can be applied to our work. First, instead of a first-order model, one can choose to apply a second-order model that maps headway to the acceleration of the vehicle. Such a model would be more physically accurate but more computationally demanding and less analytically tractable. We note that the lane-changing mechanism proposed in this work can be directly applied to second-order models and expect future work to be done in this direction.

Second, the current choices of $r$ and $p$ that affect how fast the frustration level changes can be better calibrated with more experimental results. For example, surveys of drivers may be used to examine how being passed or being in a slower lane affects their lane-changing decision. Additionally, functional magnetic resonance imaging (fMRI) may be used to monitor brain activities during or after frustration~\cite{linke2023persistent}. Similar techniques may be used on participants using driving simulations in laboratory conditions.

Third, it would be interesting to see how our proposed lane-changing mechanism responds to lanes of different characteristics such as fast lanes and slow lanes. 
One way to introduce such asymmetry could be to use different values of maximal velocity $V$ in different lanes. 
Modeling mandatory lane changes is another direction of extension. It is straightforward to implement ideas similar to \citet{yang1996microscopic} that model the probability of lane change as a function of the distance to the downstream node in our framework.

It is also worth pointing out that our results of the minimal velocity gain of the underlying car-following model may only applied to roads of relatively uniform density, which is a common feature of urban traffic. We suspect that the velocity gain of aggressive driver could be more pronounced in the roads where vehicles are moving in isolated ``packs'', commonly seen in long-distance roads connecting major cities. While the origin of the pack formation is unclear, one way to produce the scenario in a numerical simulation is to introduce vehicles whose maximal velocities are small. It would be interesting to see how much better an aggressive driver might gain under these situations.

\section*{Acknowledgement}
The authors are grateful for the generous support of Ohio Wesleyan University's Summer Science Research Program (SSRP).


%


 \bibliographystyle{elsarticle-num-names} 
 \bibliography{ref.bib}

\begin{thebibliography}{63}
\expandafter\ifx\csname natexlab\endcsname\relax\def\natexlab#1{#1}\fi
\providecommand{\url}[1]{\texttt{#1}}
\providecommand{\href}[2]{#2}
\providecommand{\path}[1]{#1}
\providecommand{\DOIprefix}{doi:}
\providecommand{\ArXivprefix}{arXiv:}
\providecommand{\URLprefix}{URL: }
\providecommand{\Pubmedprefix}{pmid:}
\providecommand{\doi}[1]{\href{http://dx.doi.org/#1}{\path{#1}}}
\providecommand{\Pubmed}[1]{\href{pmid:#1}{\path{#1}}}
\providecommand{\bibinfo}[2]{#2}
\ifx\xfnm\relax \def\xfnm[#1]{\unskip,\space#1}\fi
\bibitem[{Greenshields et~al.(1935)Greenshields, Bibbins, Channing, and
  Miller}]{greenshields1935study}
\bibinfo{author}{B.~D. Greenshields}, \bibinfo{author}{J.~R. Bibbins},
  \bibinfo{author}{W.~Channing}, \bibinfo{author}{H.~H. Miller},
\newblock \bibinfo{title}{A study of traffic capacity},
\newblock in: \bibinfo{booktitle}{Highway research board proceedings},
  volume~\bibinfo{volume}{14}, \bibinfo{organization}{Washington, DC},
  \bibinfo{year}{1935}, pp. \bibinfo{pages}{448--477}.
\bibitem[{Lighthill and Whitham(1955)}]{lighthill1955kinematic}
\bibinfo{author}{M.~J. Lighthill}, \bibinfo{author}{G.~B. Whitham},
\newblock \bibinfo{title}{On kinematic waves ii. a theory of traffic flow on
  long crowded roads},
\newblock \bibinfo{journal}{Proceedings of the royal society of london. series
  a. mathematical and physical sciences} \bibinfo{volume}{229}
  (\bibinfo{year}{1955}) \bibinfo{pages}{317--345}.
\bibitem[{Richards(1956)}]{richards1956shock}
\bibinfo{author}{P.~I. Richards},
\newblock \bibinfo{title}{Shock waves on the highway},
\newblock \bibinfo{journal}{Operations research} \bibinfo{volume}{4}
  (\bibinfo{year}{1956}) \bibinfo{pages}{42--51}.
\bibitem[{Chandler et~al.(1958)Chandler, Herman, and
  Montroll}]{chandler1958traffic}
\bibinfo{author}{R.~E. Chandler}, \bibinfo{author}{R.~Herman},
  \bibinfo{author}{E.~W. Montroll},
\newblock \bibinfo{title}{Traffic dynamics: studies in car following},
\newblock \bibinfo{journal}{Operations research} \bibinfo{volume}{6}
  (\bibinfo{year}{1958}) \bibinfo{pages}{165--184}.
\bibitem[{Helbing(2001)}]{helbing2001traffic}
\bibinfo{author}{D.~Helbing},
\newblock \bibinfo{title}{Traffic and related self-driven many-particle
  systems},
\newblock \bibinfo{journal}{Reviews of modern physics} \bibinfo{volume}{73}
  (\bibinfo{year}{2001}) \bibinfo{pages}{1067}.
\bibitem[{Treiber and Kesting(2013)}]{treiber2013traffic}
\bibinfo{author}{M.~Treiber}, \bibinfo{author}{A.~Kesting},
\newblock \bibinfo{title}{Traffic flow dynamics},
\newblock \bibinfo{journal}{Traffic Flow Dynamics: Data, Models and Simulation,
  Springer-Verlag Berlin Heidelberg} \bibinfo{volume}{227}
  (\bibinfo{year}{2013}) \bibinfo{pages}{228}.
\bibitem[{Levinson and Krizek(2007)}]{levinson2007planning}
\bibinfo{author}{D.~M. Levinson}, \bibinfo{author}{K.~J. Krizek},
  \bibinfo{title}{Planning for place and plexus: Metropolitan land use and
  transport}, \bibinfo{publisher}{Routledge}, \bibinfo{year}{2007}.
\bibitem[{Nagel and Schreckenberg(1992)}]{nagel1992cellular}
\bibinfo{author}{K.~Nagel}, \bibinfo{author}{M.~Schreckenberg},
\newblock \bibinfo{title}{A cellular automaton model for freeway traffic},
\newblock \bibinfo{journal}{Journal de physique I} \bibinfo{volume}{2}
  (\bibinfo{year}{1992}) \bibinfo{pages}{2221--2229}.
\bibitem[{Toledo(2007)}]{toledo2007driving}
\bibinfo{author}{T.~Toledo},
\newblock \bibinfo{title}{Driving behaviour: models and challenges},
\newblock \bibinfo{journal}{Transport Reviews} \bibinfo{volume}{27}
  (\bibinfo{year}{2007}) \bibinfo{pages}{65--84}.
\bibitem[{Saberi et~al.(2020)Saberi, Hamedmoghadam, Ashfaq, Hosseini, Gu,
  Shafiei, Nair, Dixit, Gardner, Waller et~al.}]{saberi2020simple}
\bibinfo{author}{M.~Saberi}, \bibinfo{author}{H.~Hamedmoghadam},
  \bibinfo{author}{M.~Ashfaq}, \bibinfo{author}{S.~A. Hosseini},
  \bibinfo{author}{Z.~Gu}, \bibinfo{author}{S.~Shafiei}, \bibinfo{author}{D.~J.
  Nair}, \bibinfo{author}{V.~Dixit}, \bibinfo{author}{L.~Gardner},
  \bibinfo{author}{S.~T. Waller}, et~al.,
\newblock \bibinfo{title}{A simple contagion process describes spreading of
  traffic jams in urban networks},
\newblock \bibinfo{journal}{Nature communications} \bibinfo{volume}{11}
  (\bibinfo{year}{2020}) \bibinfo{pages}{1616}.
\bibitem[{Amb{\"u}hl et~al.(2023)Amb{\"u}hl, Menendez, and
  Gonz{\'a}lez}]{ambuhl2023understanding}
\bibinfo{author}{L.~Amb{\"u}hl}, \bibinfo{author}{M.~Menendez},
  \bibinfo{author}{M.~C. Gonz{\'a}lez},
\newblock \bibinfo{title}{Understanding congestion propagation by combining
  percolation theory with the macroscopic fundamental diagram},
\newblock \bibinfo{journal}{Communications Physics} \bibinfo{volume}{6}
  (\bibinfo{year}{2023}) \bibinfo{pages}{26}.
\bibitem[{Kometani(1959)}]{kometani1959dynamic}
\bibinfo{author}{E.~Kometani},
\newblock \bibinfo{title}{Dynamic behavior of traffic with a nonlinear
  spacing-speed relationship},
\newblock \bibinfo{journal}{Theory of Traffic Flow (Proc. of Sym. on TTF (GM))}
   (\bibinfo{year}{1959}) \bibinfo{pages}{105--119}.
\bibitem[{Newell(1961)}]{newell1961nonlinear}
\bibinfo{author}{G.~F. Newell},
\newblock \bibinfo{title}{Nonlinear effects in the dynamics of car following},
\newblock \bibinfo{journal}{Operations research} \bibinfo{volume}{9}
  (\bibinfo{year}{1961}) \bibinfo{pages}{209--229}.
\bibitem[{Treiterer and Myers(1974)}]{treiterer1974hysteresis}
\bibinfo{author}{J.~Treiterer}, \bibinfo{author}{J.~Myers},
\newblock \bibinfo{title}{The hysteresis phenomenon in traffic flow},
\newblock \bibinfo{journal}{Transportation and traffic theory}
  \bibinfo{volume}{6} (\bibinfo{year}{1974}) \bibinfo{pages}{13--38}.
\bibitem[{Bekey et~al.(1977)Bekey, Burnham, and Seo}]{bekey1977control}
\bibinfo{author}{G.~A. Bekey}, \bibinfo{author}{G.~O. Burnham},
  \bibinfo{author}{J.~Seo},
\newblock \bibinfo{title}{Control theoretic models of human drivers in car
  following},
\newblock \bibinfo{journal}{Human Factors} \bibinfo{volume}{19}
  (\bibinfo{year}{1977}) \bibinfo{pages}{399--413}.
\bibitem[{Aron(1988)}]{aron1988car}
\bibinfo{author}{M.~Aron},
\newblock \bibinfo{title}{Car following in an urban network: simulation and
  experiments},
\newblock \bibinfo{journal}{Planning and Transport Research and Computation}
  (\bibinfo{year}{1988}).
\bibitem[{Gipps(1981)}]{gipps1981behavioural}
\bibinfo{author}{P.~G. Gipps},
\newblock \bibinfo{title}{A behavioural car-following model for computer
  simulation},
\newblock \bibinfo{journal}{Transportation research part B: methodological}
  \bibinfo{volume}{15} (\bibinfo{year}{1981}) \bibinfo{pages}{105--111}.
\bibitem[{Gipps(1986)}]{gipps1986model}
\bibinfo{author}{P.~G. Gipps},
\newblock \bibinfo{title}{A model for the structure of lane-changing
  decisions},
\newblock \bibinfo{journal}{Transportation Research Part B: Methodological}
  \bibinfo{volume}{20} (\bibinfo{year}{1986}) \bibinfo{pages}{403--414}.
\bibitem[{Ozaki(1993)}]{ozaki1993reaction}
\bibinfo{author}{H.~Ozaki},
\newblock \bibinfo{title}{Reaction and anticipation in the car following
  behavior},
\newblock in: \bibinfo{booktitle}{Proceedings of the 13-th International
  Symposium on Traffic and Transportation Theory, 1993}, \bibinfo{year}{1993}.
\bibitem[{Yang and Koutsopoulos(1996)}]{yang1996microscopic}
\bibinfo{author}{Q.~I. Yang}, \bibinfo{author}{H.~N. Koutsopoulos},
\newblock \bibinfo{title}{A microscopic traffic simulator for evaluation of
  dynamic traffic management systems},
\newblock \bibinfo{journal}{Transportation Research Part C: Emerging
  Technologies} \bibinfo{volume}{4} (\bibinfo{year}{1996})
  \bibinfo{pages}{113--129}.
\bibitem[{Chen et~al.(2012)Chen, Laval, Zheng, and Ahn}]{chen2012behavioral}
\bibinfo{author}{D.~Chen}, \bibinfo{author}{J.~Laval},
  \bibinfo{author}{Z.~Zheng}, \bibinfo{author}{S.~Ahn},
\newblock \bibinfo{title}{A behavioral car-following model that captures
  traffic oscillations},
\newblock \bibinfo{journal}{Transportation research part B: methodological}
  \bibinfo{volume}{46} (\bibinfo{year}{2012}) \bibinfo{pages}{744--761}.
\bibitem[{Kerner and Rehborn(1996)}]{kerner1996experimental}
\bibinfo{author}{B.~S. Kerner}, \bibinfo{author}{H.~Rehborn},
\newblock \bibinfo{title}{Experimental properties of complexity in traffic
  flow},
\newblock \bibinfo{journal}{Physical Review E} \bibinfo{volume}{53}
  (\bibinfo{year}{1996}) \bibinfo{pages}{R4275}.
\bibitem[{Kerner(1998)}]{kerner1998experimental}
\bibinfo{author}{B.~S. Kerner},
\newblock \bibinfo{title}{Experimental features of self-organization in traffic
  flow},
\newblock \bibinfo{journal}{Physical review letters} \bibinfo{volume}{81}
  (\bibinfo{year}{1998}) \bibinfo{pages}{3797}.
\bibitem[{Bando et~al.(1994)Bando, Hasebe, Nakayama, Shibata, and
  Sugiyama}]{bando1994structure}
\bibinfo{author}{M.~Bando}, \bibinfo{author}{K.~Hasebe},
  \bibinfo{author}{A.~Nakayama}, \bibinfo{author}{A.~Shibata},
  \bibinfo{author}{Y.~Sugiyama},
\newblock \bibinfo{title}{Structure stability of congestion in traffic
  dynamics},
\newblock \bibinfo{journal}{Japan Journal of Industrial and Applied
  Mathematics} \bibinfo{volume}{11} (\bibinfo{year}{1994})
  \bibinfo{pages}{203--223}.
\bibitem[{Bando et~al.(1995)Bando, Hasebe, Nakayama, Shibata, and
  Sugiyama}]{bando1995dynamical}
\bibinfo{author}{M.~Bando}, \bibinfo{author}{K.~Hasebe},
  \bibinfo{author}{A.~Nakayama}, \bibinfo{author}{A.~Shibata},
  \bibinfo{author}{Y.~Sugiyama},
\newblock \bibinfo{title}{Dynamical model of traffic congestion and numerical
  simulation},
\newblock \bibinfo{journal}{Physical review E} \bibinfo{volume}{51}
  (\bibinfo{year}{1995}) \bibinfo{pages}{1035}.
\bibitem[{Holland(1998)}]{holland1998generalised}
\bibinfo{author}{E.~Holland},
\newblock \bibinfo{title}{A generalised stability criterion for motorway
  traffic},
\newblock \bibinfo{journal}{Transportation Research Part B: Methodological}
  \bibinfo{volume}{32} (\bibinfo{year}{1998}) \bibinfo{pages}{141--154}.
\bibitem[{Kerner(2009)}]{kerner2009introduction}
\bibinfo{author}{B.~S. Kerner}, \bibinfo{title}{Introduction to modern traffic
  flow theory and control: the long road to three-phase traffic theory},
  \bibinfo{publisher}{Springer Science \& Business Media},
  \bibinfo{year}{2009}.
\bibitem[{Wilson(2008)}]{wilson2008mechanisms}
\bibinfo{author}{R.~E. Wilson},
\newblock \bibinfo{title}{Mechanisms for spatio-temporal pattern formation in
  highway traffic models},
\newblock \bibinfo{journal}{Philosophical Transactions of the Royal Society A:
  Mathematical, Physical and Engineering Sciences} \bibinfo{volume}{366}
  (\bibinfo{year}{2008}) \bibinfo{pages}{2017--2032}.
\bibitem[{Seibold et~al.(2012)Seibold, Flynn, Kasimov, and
  Rosales}]{seibold2012constructing}
\bibinfo{author}{B.~Seibold}, \bibinfo{author}{M.~R. Flynn},
  \bibinfo{author}{A.~R. Kasimov}, \bibinfo{author}{R.~R. Rosales},
\newblock \bibinfo{title}{Constructing set-valued fundamental diagrams from
  jamiton solutions in second order traffic models},
\newblock \bibinfo{journal}{arXiv preprint arXiv:1204.5510}
  (\bibinfo{year}{2012}).
\bibitem[{Newell(1962)}]{newell1962theories}
\bibinfo{author}{G.~F. Newell},
\newblock \bibinfo{title}{Theories of instability in dense highway traffic},
\newblock \bibinfo{journal}{J. Operations Research Society of Japan}
  \bibinfo{volume}{5} (\bibinfo{year}{1962}) \bibinfo{pages}{9--54}.
\bibitem[{Kesting and Treiber(2008)}]{kesting2008reaction}
\bibinfo{author}{A.~Kesting}, \bibinfo{author}{M.~Treiber},
\newblock \bibinfo{title}{How reaction time, update time, and adaptation time
  influence the stability of traffic flow},
\newblock \bibinfo{journal}{Computer-Aided Civil and Infrastructure
  Engineering} \bibinfo{volume}{23} (\bibinfo{year}{2008})
  \bibinfo{pages}{125--137}.
\bibitem[{Laval and Daganzo(2006)}]{laval2006lane}
\bibinfo{author}{J.~A. Laval}, \bibinfo{author}{C.~F. Daganzo},
\newblock \bibinfo{title}{Lane-changing in traffic streams},
\newblock \bibinfo{journal}{Transportation Research Part B: Methodological}
  \bibinfo{volume}{40} (\bibinfo{year}{2006}) \bibinfo{pages}{251--264}.
\bibitem[{Laval et~al.(2014)Laval, Toth, and Zhou}]{laval2014parsimonious}
\bibinfo{author}{J.~A. Laval}, \bibinfo{author}{C.~S. Toth},
  \bibinfo{author}{Y.~Zhou},
\newblock \bibinfo{title}{A parsimonious model for the formation of
  oscillations in car-following models},
\newblock \bibinfo{journal}{Transportation Research Part B: Methodological}
  \bibinfo{volume}{70} (\bibinfo{year}{2014}) \bibinfo{pages}{228--238}.
\bibitem[{Shang et~al.(2020)Shang, Hauer, and Stern}]{shang2020cut}
\bibinfo{author}{M.~Shang}, \bibinfo{author}{F.~Hauer},
  \bibinfo{author}{R.~Stern},
\newblock \bibinfo{title}{Do cut-ins matter: Assessing the impact of lane
  changing and string stability on traffic flow},
\newblock in: \bibinfo{booktitle}{2020 IEEE 23rd International Conference on
  Intelligent Transportation Systems (ITSC)}, \bibinfo{organization}{IEEE},
  \bibinfo{year}{2020}, pp. \bibinfo{pages}{1--6}.
\bibitem[{Chauhan et~al.(2022)Chauhan, Kanagaraj, and
  Asaithambi}]{chauhan2022understanding}
\bibinfo{author}{P.~Chauhan}, \bibinfo{author}{V.~Kanagaraj},
  \bibinfo{author}{G.~Asaithambi},
\newblock \bibinfo{title}{Understanding the mechanism of lane changing process
  and dynamics using microscopic traffic data},
\newblock \bibinfo{journal}{Physica A: Statistical Mechanics and its
  Applications} \bibinfo{volume}{593} (\bibinfo{year}{2022})
  \bibinfo{pages}{126981}.
\bibitem[{Gao and Levinson(2023)}]{gao2023lane}
\bibinfo{author}{Y.~Gao}, \bibinfo{author}{D.~Levinson},
\newblock \bibinfo{title}{Lane changing and congestion are mutually
  reinforcing?},
\newblock \bibinfo{journal}{Communications in Transportation Research}
  \bibinfo{volume}{3} (\bibinfo{year}{2023}) \bibinfo{pages}{100101}.
\bibitem[{Zhang et~al.(2024)Zhang, Rong, Gao, and Chen}]{zhang2024effects}
\bibinfo{author}{K.~Zhang}, \bibinfo{author}{J.~Rong},
  \bibinfo{author}{Y.~Gao}, \bibinfo{author}{Y.~Chen},
\newblock \bibinfo{title}{Effects of lane imbalance on capacity drop and
  emission in expressway merging areas: A simulation analysis},
\newblock \bibinfo{journal}{Sustainability} \bibinfo{volume}{16}
  (\bibinfo{year}{2024}) \bibinfo{pages}{10388}.
\bibitem[{Hidas(2002)}]{hidas2002modelling}
\bibinfo{author}{P.~Hidas},
\newblock \bibinfo{title}{Modelling lane changing and merging in microscopic
  traffic simulation},
\newblock \bibinfo{journal}{Transportation Research Part C: Emerging
  Technologies} \bibinfo{volume}{10} (\bibinfo{year}{2002})
  \bibinfo{pages}{351--371}.
\bibitem[{Hidas(2005)}]{hidas2005modelling}
\bibinfo{author}{P.~Hidas},
\newblock \bibinfo{title}{Modelling vehicle interactions in microscopic
  simulation of merging and weaving},
\newblock \bibinfo{journal}{Transportation Research Part C: Emerging
  Technologies} \bibinfo{volume}{13} (\bibinfo{year}{2005})
  \bibinfo{pages}{37--62}.
\bibitem[{Kesting et~al.(2007)Kesting, Treiber, and
  Helbing}]{kesting2007general}
\bibinfo{author}{A.~Kesting}, \bibinfo{author}{M.~Treiber},
  \bibinfo{author}{D.~Helbing},
\newblock \bibinfo{title}{General lane-changing model mobil for car-following
  models},
\newblock \bibinfo{journal}{Transportation Research Record}
  \bibinfo{volume}{1999} (\bibinfo{year}{2007}) \bibinfo{pages}{86--94}.
\bibitem[{Ahmed et~al.(1996)Ahmed, Ben-Akiva, Koutsopoulos, and
  Mishalani}]{ahmed1996models}
\bibinfo{author}{K.~Ahmed}, \bibinfo{author}{M.~Ben-Akiva},
  \bibinfo{author}{H.~Koutsopoulos}, \bibinfo{author}{R.~Mishalani},
\newblock \bibinfo{title}{Models of freeway lane changing and gap acceptance
  behavior},
\newblock \bibinfo{journal}{Transportation and traffic theory}
  \bibinfo{volume}{13} (\bibinfo{year}{1996}) \bibinfo{pages}{501--515}.
\bibitem[{Ahmed(1999)}]{ahmed1999modeling}
\bibinfo{author}{K.~I. Ahmed}, \bibinfo{title}{Modeling drivers' acceleration
  and lane changing behavior}, Ph.D. thesis, Massachusetts Institute of
  Technology, \bibinfo{year}{1999}.
\bibitem[{Toledo et~al.(2003)Toledo, Koutsopoulos, and
  Ben-Akiva}]{toledo2003modeling}
\bibinfo{author}{T.~Toledo}, \bibinfo{author}{H.~N. Koutsopoulos},
  \bibinfo{author}{M.~E. Ben-Akiva},
\newblock \bibinfo{title}{Modeling integrated lane-changing behavior},
\newblock \bibinfo{journal}{Transportation Research Record}
  \bibinfo{volume}{1857} (\bibinfo{year}{2003}) \bibinfo{pages}{30--38}.
\bibitem[{Sun and Elefteriadou(2011)}]{sun2011lane}
\bibinfo{author}{D.~J. Sun}, \bibinfo{author}{L.~Elefteriadou},
\newblock \bibinfo{title}{Lane-changing behavior on urban streets: A focus
  group-based study},
\newblock \bibinfo{journal}{Applied ergonomics} \bibinfo{volume}{42}
  (\bibinfo{year}{2011}) \bibinfo{pages}{682--691}.
\bibitem[{Sun and Elefteriadou(2012)}]{sun2012lane}
\bibinfo{author}{D.~Sun}, \bibinfo{author}{L.~Elefteriadou},
\newblock \bibinfo{title}{Lane-changing behavior on urban streets: An
  “in-vehicle” field experiment-based study},
\newblock \bibinfo{journal}{Computer-Aided Civil and Infrastructure
  Engineering} \bibinfo{volume}{27} (\bibinfo{year}{2012})
  \bibinfo{pages}{525--542}.
\bibitem[{Nagatani(1993)}]{nagatani1993self}
\bibinfo{author}{T.~Nagatani},
\newblock \bibinfo{title}{Self-organization and phase transition in
  traffic-flow model of a two-lane roadway},
\newblock \bibinfo{journal}{Journal of Physics A: Mathematical and General}
  \bibinfo{volume}{26} (\bibinfo{year}{1993}) \bibinfo{pages}{L781}.
\bibitem[{Rickert et~al.(1996)Rickert, Nagel, Schreckenberg, and
  Latour}]{rickert1996two}
\bibinfo{author}{M.~Rickert}, \bibinfo{author}{K.~Nagel},
  \bibinfo{author}{M.~Schreckenberg}, \bibinfo{author}{A.~Latour},
\newblock \bibinfo{title}{Two lane traffic simulations using cellular
  automata},
\newblock \bibinfo{journal}{Physica A: Statistical Mechanics and its
  Applications} \bibinfo{volume}{231} (\bibinfo{year}{1996})
  \bibinfo{pages}{534--550}.
\bibitem[{Wagner et~al.(1997)Wagner, Nagel, and Wolf}]{wagner1997realistic}
\bibinfo{author}{P.~Wagner}, \bibinfo{author}{K.~Nagel}, \bibinfo{author}{D.~E.
  Wolf},
\newblock \bibinfo{title}{Realistic multi-lane traffic rules for cellular
  automata},
\newblock \bibinfo{journal}{Physica A: Statistical Mechanics and its
  Applications} \bibinfo{volume}{234} (\bibinfo{year}{1997})
  \bibinfo{pages}{687--698}.
\bibitem[{Nagel et~al.(1998)Nagel, Wolf, Wagner, and Simon}]{nagel1998two}
\bibinfo{author}{K.~Nagel}, \bibinfo{author}{D.~E. Wolf},
  \bibinfo{author}{P.~Wagner}, \bibinfo{author}{P.~Simon},
\newblock \bibinfo{title}{Two-lane traffic rules for cellular automata: A
  systematic approach},
\newblock \bibinfo{journal}{Physical Review E} \bibinfo{volume}{58}
  (\bibinfo{year}{1998}) \bibinfo{pages}{1425}.
\bibitem[{Maerivoet and De~Moor(2005)}]{maerivoet2005cellular}
\bibinfo{author}{S.~Maerivoet}, \bibinfo{author}{B.~De~Moor},
\newblock \bibinfo{title}{Cellular automata models of road traffic},
\newblock \bibinfo{journal}{Physics reports} \bibinfo{volume}{419}
  (\bibinfo{year}{2005}) \bibinfo{pages}{1--64}.
\bibitem[{Worrall et~al.(1970)Worrall, Bullen, and Gur}]{worrall1970elementary}
\bibinfo{author}{R.~Worrall}, \bibinfo{author}{A.~Bullen},
  \bibinfo{author}{Y.~Gur},
\newblock \bibinfo{title}{An elementary stochastic model of lane-changing on a
  multilane highway},
\newblock \bibinfo{journal}{Highway Research Record}  (\bibinfo{year}{1970}).
\bibitem[{Pentland and Liu(1999)}]{pentland1999modeling}
\bibinfo{author}{A.~Pentland}, \bibinfo{author}{A.~Liu},
\newblock \bibinfo{title}{Modeling and prediction of human behavior},
\newblock \bibinfo{journal}{Neural computation} \bibinfo{volume}{11}
  (\bibinfo{year}{1999}) \bibinfo{pages}{229--242}.
\bibitem[{Sheu and Ritchie(2001)}]{sheu2001stochastic}
\bibinfo{author}{J.-B. Sheu}, \bibinfo{author}{S.~G. Ritchie},
\newblock \bibinfo{title}{Stochastic modeling and real-time prediction of
  vehicular lane-changing behavior},
\newblock \bibinfo{journal}{Transportation Research Part B: Methodological}
  \bibinfo{volume}{35} (\bibinfo{year}{2001}) \bibinfo{pages}{695--716}.
\bibitem[{Toledo and Katz(2009)}]{toledo2009state}
\bibinfo{author}{T.~Toledo}, \bibinfo{author}{R.~Katz},
\newblock \bibinfo{title}{State dependence in lane-changing models},
\newblock \bibinfo{journal}{Transportation research record}
  \bibinfo{volume}{2124} (\bibinfo{year}{2009}) \bibinfo{pages}{81--88}.
\bibitem[{Singh and Li(2011)}]{singh2011estimation}
\bibinfo{author}{K.~Singh}, \bibinfo{author}{B.~Li},
\newblock \bibinfo{title}{Estimation of traffic densities for multilane
  roadways using a markov model approach},
\newblock \bibinfo{journal}{IEEE Transactions on industrial electronics}
  \bibinfo{volume}{59} (\bibinfo{year}{2011}) \bibinfo{pages}{4369--4376}.
\bibitem[{Zheng(2014)}]{zheng2014recent}
\bibinfo{author}{Z.~Zheng},
\newblock \bibinfo{title}{Recent developments and research needs in modeling
  lane changing},
\newblock \bibinfo{journal}{Transportation research part B: methodological}
  \bibinfo{volume}{60} (\bibinfo{year}{2014}) \bibinfo{pages}{16--32}.
\bibitem[{Smith(2011)}]{smith2011introduction}
\bibinfo{author}{H.~L. Smith}, \bibinfo{title}{An introduction to delay
  differential equations with applications to the life sciences},
  volume~\bibinfo{volume}{57}, \bibinfo{publisher}{springer New York},
  \bibinfo{year}{2011}.
\bibitem[{Kerner(1999)}]{kerner1999physics}
\bibinfo{author}{B.~S. Kerner},
\newblock \bibinfo{title}{The physics of traffic},
\newblock \bibinfo{journal}{Physics World} \bibinfo{volume}{12}
  (\bibinfo{year}{1999}) \bibinfo{pages}{25}.
\bibitem[{Chen et~al.(2014)Chen, Ahn, Laval, and Zheng}]{chen2014periodicity}
\bibinfo{author}{D.~Chen}, \bibinfo{author}{S.~Ahn},
  \bibinfo{author}{J.~Laval}, \bibinfo{author}{Z.~Zheng},
\newblock \bibinfo{title}{On the periodicity of traffic oscillations and
  capacity drop: The role of driver characteristics},
\newblock \bibinfo{journal}{Transportation research part B: methodological}
  \bibinfo{volume}{59} (\bibinfo{year}{2014}) \bibinfo{pages}{117--136}.
\bibitem[{James(2009)}]{james2009road}
\bibinfo{author}{L.~James}, \bibinfo{title}{Road rage and aggressive driving:
  Steering clear of highway warfare}, \bibinfo{publisher}{Prometheus Books},
  \bibinfo{year}{2009}.
\bibitem[{Zheng et~al.(2010)Zheng, Ahn, and Monsere}]{zheng2010impact}
\bibinfo{author}{Z.~Zheng}, \bibinfo{author}{S.~Ahn}, \bibinfo{author}{C.~M.
  Monsere},
\newblock \bibinfo{title}{Impact of traffic oscillations on freeway crash
  occurrences},
\newblock \bibinfo{journal}{Accident Analysis \& Prevention}
  \bibinfo{volume}{42} (\bibinfo{year}{2010}) \bibinfo{pages}{626--636}.
\bibitem[{Zheng et~al.(2013)Zheng, Ahn, Chen, and Laval}]{zheng2013effects}
\bibinfo{author}{Z.~Zheng}, \bibinfo{author}{S.~Ahn},
  \bibinfo{author}{D.~Chen}, \bibinfo{author}{J.~Laval},
\newblock \bibinfo{title}{The effects of lane-changing on the immediate
  follower: Anticipation, relaxation, and change in driver characteristics},
\newblock \bibinfo{journal}{Transportation research part C: emerging
  technologies} \bibinfo{volume}{26} (\bibinfo{year}{2013})
  \bibinfo{pages}{367--379}.
\bibitem[{Linke et~al.(2023)Linke, Haller, Xu, Nguyen, Chue, Botz-Zapp,
  Revzina, Perlstein, Ross, Tseng et~al.}]{linke2023persistent}
\bibinfo{author}{J.~O. Linke}, \bibinfo{author}{S.~P. Haller},
  \bibinfo{author}{E.~P. Xu}, \bibinfo{author}{L.~T. Nguyen},
  \bibinfo{author}{A.~E. Chue}, \bibinfo{author}{C.~Botz-Zapp},
  \bibinfo{author}{O.~Revzina}, \bibinfo{author}{S.~Perlstein},
  \bibinfo{author}{A.~J. Ross}, \bibinfo{author}{W.-L. Tseng}, et~al.,
\newblock \bibinfo{title}{Persistent frustration-induced reconfigurations of
  brain networks predict individual differences in irritability},
\newblock \bibinfo{journal}{Journal of the American Academy of Child \&
  Adolescent Psychiatry} \bibinfo{volume}{62} (\bibinfo{year}{2023})
  \bibinfo{pages}{684--695}.

\end{thebibliography}






\end{document}